


\baselineskip=14pt
\parskip=10pt

\magnification=\magstephalf

\def\S{{\cal S}}
\def\1{{\overline{1}}}
\def\2{{\overline{2}}}
\parindent=0pt
\overfullrule=0in

\def\frac#1#2{{#1 \over #2}}
\centerline
{\bf  
Counting (and Randomly Generating) Hamiltonian Cycles in Rectangular Grids
}
\bigskip
\centerline
{\it Pablo BLANCO and Doron ZEILBERGER}

\bigskip

\qquad\qquad  {\it Dedicated to Volker Strehl (born June 23, 1945), great enumerator and dear friend, on his (forthcoming) ``perfect fourth power birthday"}

{\bf Abstract}:  We first fully implement, in Maple, the ingenious method of Robert Stoyan and Volker Strehl from 1995 to automatically derive generating functions for the number of
Hamiltonian cycles in an $m$ by $n$ grid graph, for a fixed width $m$, but general length $n$, and actually compute these generating functions for all $ 2 \leq m \leq 10$.
We also show how to  generate a uniformly-at-random such Hamiltonian cycle, and give a full implementation in Maple.
Finally we derive more informative generating functions for other parameters besides the length of the underlying grid graph.

{\bf Preface: How it all started}

In Nov. 23, 2025 the New Your Times Sunday magazine started posing a new kind of puzzle that they called {\it Spell-Weaving}.
It was invented by puzzle-making whiz Rudolfo Kurchan, and discontinued Feb. 22, 2026 (replaced by another kind of puzzle called {\it Recto}).

You are given a $5\times5$ grid (``chessboard'') (later weeks had
$6 \times 6$ grids and the final week had a $7 \times 7$ grid) as well as a word. The letters of the word are scattered in the grid, and the solver is asked to
{\it spell that word}, by starting at the first letter, ending at the last letter, and {\it passing through every square exactly once}.
In other words, find a {\it Hamiltonian path} in the given grid graph that passes  through the designated letters in the right order.

Here is the solution to the very challenging (last one!) Feb. 22, 2026 puzzle:

{\tt https://sites.math.rutgers.edu/\~{}zeilberg/tokhniot/oKurchanComplicated.png} \quad .

This led us to write a Maple package:

{\tt https://sites.math.rutgers.edu/\~{}zeilberg/tokhniot/Kurchan.txt} \quad,

that not only solves these kinds of puzzles, but also creates them, leading to the following web puzzle book:

{\tt https://sites.math.rutgers.edu/\~{}zeilberg/SpellWeaving/SpellWeaving.html} \quad .

{\bf How Many?}

Being enumerators, we asked ``how many Hamiltonian cycles are there in a given grid graph?''. Of course, we are not naive, and we know, {\it for sure}, that
we will {\bf never know} the exact number of Hamiltonian cycles on the graph $P_{100} \times P_{100}$. But how about the exact number of
Hamiltonian cycles on the grid graph $P_{10} \times P_{100}$? After reading this article, the reader will find out how to obtain it quickly.
It turns out that this number is {\bf exactly}

$$
\matrix
{
2841755307998403180696485173480879907420461708673514070665759422586711\backslash\cr
26855416799214435461577164935511762299757966788827828321166383429987198
} \quad .
$$

For the corresponding numbers for $P_{10} \times P_{1000}$  and $P_{10} \times P_{10000}$  see here:

{\tt http://sites.math.rutgers.edu/\~{}zeilberg/tokhniot/oVolker7.txt} \quad .

The latter is approximately: $8.399066204805426684770915677726152158842\cdot 10^{14310}$. Way way more than the number of atoms in the universe!

We first started modestly, and using {\it dynamical programming}, cranked out the number of Hamiltonian cycles for $P_3 \times P_n$, for $1 \leq n \leq 6$
(these turned out to be powers of $2$ and a nice exercise to prove),
and then we computed these numbers for  $P_4 \times P_n$, for $1 \leq n \leq 6$. {\it Of course}, we went to our favorite website, the OEIS [Sl], and found out that people
already worked on it before:

{\tt https://oeis.org/A006864}.

That's how we came across references [F] and [Kw]. Later, also thanks to the OEIS, we came across the paper [SS] and we were really impressed!
It was  a nice surprise that one of the authors, Volker Strehl, is a an old friend of one of us (DZ), but he never told us about that paper!
It was co-authored with Robert Stoyan, in what seems to be  a  published version of a {\it diplomarbeit}, written by Stoyan and supervised by Strehl.

{\bf Implementing the Stoyan-Strehl Approach}

We refer the reader to [SS] for a description of the elegant Stoyan-Strehl approach. Unfortunately, their article was not accompanied by the
supporting computer code that was used  to implement their method and to output the explicit generating functions found at the end
of their paper (they gave the cases $2\leq m \leq 6$ completely, and parts of the cases $7 \leq m \leq 8$, listing the degrees of the
numerators and denominators of the enumerating generating functions). Since we couldn't find an extant implementation on the internet,
our first goal was to fill this regrettable {\it gap}, and have a full, user-friendly, freely available implementation, in our favorite
computer algebra system, {\tt Maple}. So we wrote the Maple package

{\tt https://sites.math.rutgers.edu/\~{}zeilberg/tokhniot/Volker.txt} \quad,

that fully implements the Stoyan-Strehl method, and, as we hope to show later, extend it to do much more than just crank-out the enumerating generating functions.

\vfill\eject

{\bf A brief description of the Stoyan-Strehl method}

Robert Stoyan and Volker Strehl  preferred to use the language of {\bf finite automata}. We prefer to
use the equivalent language of {\it directed graphs}. We assume that the reader read and understood the first ten pages of [SS]. Since this is
so well-written, and readily available on-line, we will only describe how we adapted their automata to the language of directed graphs.

Recall that in [SS], they establish a natural one-to-one mapping from the set of Hamiltonian cycles in $P_{m+1} \times P_{n+1}$ and 
a certain class of $m \times n$ binary matrices, obeying certain conditions.
The states of the finite automta, (or the set of vertices) are pairs $[v,\S]$ where $v$ is a member of $\{0,1\}^m$ (not all  binary vectors of length $m$ show up, only some of them), and $\S$ is some set-partition of
the subset of $\{1, \dots, m\}$ where $v_i=1$. If $v$ is a column of the matrix $M$, then the set-partition $\S$ is such that
if $v_i=1$ and $v_j=1$, ($1 \leq i < j \leq m$) $i$ and $j$ belong to the same component of $\S$ if and only if there is a path travelling from $v_i=1$ to $v_j=1$ via all the $1$-s in $M$
that lie to the left of that column. So, in the leftmost column of the good matrix $M$, only the `$1$'s that lie in the same block are `roommates', while for the rightmost column all the `$1$'s
are roommates.

To see the `alphabet' (or ``states'') of the Stoyan-Strehl language in {\tt Volker.txt}, for $m \times n$ binary matrices (corresponding to Hamiltonian cycles on $P_{m+1} \times P_{n+1}$) use the
command {\tt AlphaSS(m)}.

For example, {\tt AlphaSS(5);} would give you the alphabet consisting of $32$ `letters' of the form $[v,\S]$. The set of $v$-s that show up (out of the $2^5=32$ possible binary vectors of length $5$) happens to be:

$$
\matrix
{\{[0, 0, 0, 0, 1], [0, 0, 0, 1, 0], [0, 0, 1, 0, 0], [0, 0, 1, 0, 1], [0, 0, 1, 1, 1], [0, 1, 0, 0, 0], [0, 1, 0, 0, 1], [0, 1, 0, 1, 0], \cr
[0, 1, 1, 1, 0],  [1, 0, 0, 0, 0], [1, 0, 0, 0, 1], [1, 0, 0, 1, 0], [1, 0, 1, 0, 0], [1, 0, 1, 0, 1], [1, 0, 1, 1, 1], [1, 1, 0, 1, 1], \cr
[1, 1, 1, 0, 0], [1, 1, 1, 0, 1],[1, 1, 1, 1, 1]\} } \quad .
$$

By pure coincidence the cardinality of {\tt AlphaSS(5)} happens to be $32$. For example, there are two members $[v,\S]$ with $v=[1,1,0,1,1]$, one is
$$
[[1,1,0,1,1], \{\{1,2\}, \{4,5\}\}] \quad,
$$
where the first block of consecutive ones is {\bf not} connected to the second block of ones, and
$$
[[1,1,0,1,1], \{\{1,2,4,5\}\}] \quad,
$$
where all the ones of $[1,1,0,1,1]$ are connected to each other.

So we have the {\it alphabet}, but  what about the {\it grammar}?. Using the [SS] compatibility conditions, and the natural `connectivity', the computer can easily decide which
letters can follow any given letter, and generate the `type three grammar',  completely defining the language. This is implemented in
procedure {\tt Followers}. For example for the member 

$$
[[1,1,0,1,1], \{\{1,2\},\{4,5\}\}] \quad,
$$
of {\tt AlphaSS(5)},
$$
Followers([[1,1,0,1,1], \{\{1,2\},\{4,5\}\}]); \quad,
$$
gives
$$
\matrix{
\{[[0, 1, 0, 0, 1], \{\{2\}, \{5\}\}], & [[0, 1, 0, 1, 0], \{\{2\}, \{4\}\}], & [[0, 1, 1, 1, 0], \{\{2, 3, 4\}\}],  & [[1, 0, 0, 0, 1], \{\{1\},\{5\}\}], \cr
[[1, 0, 0, 1, 0], \{\{1\}, \{4\}\}]\}  &    &  &
}.
$$

There is a certain subset of {\tt AlphaSS(m)}, that are allowed to be the starting letters (those where only the ones in the same block are in the connected component), they are given by {\tt Starters(m)}.
For example, {\tt Starters(m)}. For example, {\tt Starters(4);} gives:
$$
\{[[1, 0, 1, 1], \{\{1\}, \{3, 4\}\}], [[1, 1, 0, 1], \{\{4\}, \{1, 2
\}\}], [[1, 1, 1, 1], \{\{1, 2, 3, 4\}\}]\} \quad .
$$

Another subset of {\tt AlphaSS(m)} is {\tt Enders(m)}, those letters that can end a Stoyan-Strehl word. For example, {\tt Enders(4);} yields

$$
\{[[1, 0, 1, 1], \{\{1, 3, 4\}\}], [[1, 1, 0, 1], \{\{1, 2, 4\}\}], [[1, 1, 1, 1], \{\{1, 2, 3, 4\}\}]\} \quad .
$$

Now that we have the alphabet and the grammar, we can model it as a directed graph.
Say that there $N_m$ such letters.
We number the letters of {\tt AlphaSS(m)} arbitrarily from $2$ to  $N_m+1$ and
add two more letters: $1$ corresponding to {\tt START}; and $N_m+2$ corresponding to {\tt END}, and translate the information in {\tt Followers} to find the corresponding edges, converting
it to a certain {\it directed graph}, with $N_m+2$ vertices. We represent it as list of sets, $G$, where $G[i]$ is the set of vertices $j$ such that there is a (directed) edge between
vertex $i$ and vertex $j$. Procedure {\tt SSdg(m)}, converts the grammar into an abstract directed graph with our convention, followed by the
`legend'. For example, {\tt SSdg(3)[1];} yields
$$
[\{6, 7\}, \{6, 7\}, \{7\}, \{6, 7\}, \{2, 4, 5, 8\}, \{6, 7\}, \{2, 3, 4, 5, 8\}, \{\}] \quad ,
$$
which denotes a directed graph with $8$ vertices, labeled $1 \dots 8$, and the outgoing neighbors of vertex $1$ are vertices $6$ and $7$, the outgoing  neighbors of vertex $2$ are also $6$ and $7$,
the only outgoing neighbors of vertex $3$ is vertex $7$, and the set of outgoing vertices of vertex $8$ is the empty set (since the last vertex is always a sink by our convention).

The second component of {\tt SSdg(m)}, {\tt SSdg(m)[2]}, gives the `dictionary', what each numbered vertex stands for in {\tt AlphaSS(m)}. For example,
{\tt SSdg(3)[2];} outputs
$$
\matrix{[1, [[0, 0, 1], \{\{3\}\}], [[0, 1, 0], \{\{2\}\}], [[1, 0, 0], \{\{1\}\}],[[1, 0, 1], \{\{1, 3\}\}], [[1, 0, 1], \{\{1\}, \{3\}\}], \cr
[[1, 1, 1], \{\{1, 2, 3\}\}], 8] \quad .}
$$

This means that the vertex labeled $1$ is still called $1$, corresponding to our convention that it is the starting vertex, {\tt START}.
Vertex $2$ stands  for the letter $[[0, 0, 1], \{\{3\}\}]$, etc.

{\bf Generating functions}

After the set-up, we use the standard {\it transfer matrix method} (e.g. [St], section 4.7), to find the generating function (or {\it weight-enumerator}), in the formal variable $x$, for
walks according to the length. As is well-known and easy to see, if the {\it adjacency matrix} of our directed graph is $M$, then the weight-enumerator of the
set of walks from vertex $i$ to vertex $j$ is the $(i,j)$ entry in the matrix.
$$
(I-xM)^{-1} \quad .
$$

For our application we are only interested in the weight-enumerator of the set of walks from the first vertex to the last one (i.e. from vertex $1$ to vertex $|AlphaSS(m)|+2$). This is implemented, for
arbitrary directed graphs, in the other Maple package accompanying this article, {\tt WDG.txt}, available from:

{\tt https://sites.math.rutgers.edu/\~{}zeilberg/tokhniot/WDG.txt} \quad .

{\bf Output}

Going back to Hamiltonian cycles in grid graphs, procedure {\tt GFnz(m,z)} outputs the rational function, in the variable $z$, whose coefficient of $z^n$ is the
exact number of Hamiltonian cycles in the grid graph $P_m \times P_n$. For example, {\tt GFnz(4,z);}, yields
$$
-\frac{z^{2}}{z^{4}-2 z^{3}+2 z^{2}+2 z -1} \quad,
$$
while {\tt GFnz(5,z);}, yields
$$
-\frac{z^{2} \left(3 z^{2}+1\right)}{2 z^{6}+11 z^{2}-1} \quad .
$$
This agrees with the output on p. 19 of [SS].

As $m$ grows larger, naturally, {\tt GFnz(m,z);} gets slower, so we have another procedure,

{\tt FGnzPC(m,z);}

that outputs the {\it pre-computed} values of these generating functions. For example,  if you type
$$
f:=GFnzPC(10,z): \quad,
$$
(make sure to have a colon, {\bf not} a semi-colon, since this rational function is huge, unless, you want to see it), followed by

{\tt taylor(f,z=0,11);}

would output
$$
z^{2}+16 z^{3}+1517 z^{4}+18684 z^{5}+1024028 z^{6}+17066492 z^{7}+681728204 z^{8}+13916993782 z^{9}+467260456608 z^{10}+ \dots  \quad,
$$
which tells you that the number of Hamiltonian cycles of the grid graph  $P_2 \times P_{10}$ is $1$, of $P_3 \times P_{10}$ is $16$,  $\dots$, and the number
of Hamiltonian cycles of the grid graph $P_{10} \times P_{10}$ is $467260456608$. To get the number of Hamiltonian cycles in the grid graph $P_{10} \times P_{100}$, type:

{\tt coeff(taylor(f,z=0,101),z,100);} \quad,

and you would get the number given at the beginning of this paper.

{\bf Weighted Counting}

But what about a generating function whose coefficient of $z^n$ is not just the total number of Hamiltonian cycles of the grid graph $P_m \times P_n$ but rather the
{\it weight-enumerator} according to some `statistic', for example, the number of edges visited at the top?, or the bottom?, or both?. The beauty of symbolic
computation is that it is so easy to ``tweak'' the code to get much more informative generating functions, that are now rational functions of not just $z$,  but also of other
variable(s), that keep(s) track of the given statistic(s). This is implemented in procedure {\tt GFnzG(m,z,w)}. If you leave $w$ as a symbol, then
you would get the rational function in $z$ and $w_1, \dots, w_{m-1}$ such that when you Taylor it at $z=0$ (i.e. {\it Maclaurin it}), and then extract the
coefficient of $z^n$ and then, in turn,  extract the coefficient of the monomial $w_1^{a_1} \cdots w_{m-1}^{a_{m-1}}$, you would get the {\bf exact} number of Hamiltonian cycles
of the grid graph $P_m \times P_n$, whose corresponding Stoyan-Strehl $(m-1) \times (n-1)$ zero-one matrices have $a_1$ ones in the first row,
$a_2$ ones in the second row, $\dots$, $a_{m-1}$ ones in the $(m-1)^{th}$ row. For example, typing (in {\tt Volker.txt})

{\tt f:=GFnzG(4,z,w);}

yields
$$
-\frac{z^{2} w_{1} w_{2} w_{3}}{z^{4} w_{1}^{3} w_{2}^{2} w_{3}^{3}-2 z^{3} w_{1}^{2} w_{2}^{2} w_{3}^{2}+z^{2} w_{1}^{2} w_{2} w_{3}-z^{2} w_{1}^{2} w_{3}^{2}+z^{2} w_{1} w_{2}^{2} w_{3}+z^{2} w_{1} w_{2} w_{3}^{2}+2 z w_{1} w_{3}-1} \quad.
$$

Now typing

{\tt expand(coeff(taylor(f,z=0,11),z,10));} \quad,

would give you the weight-enumerator of the set of Hamiltonian cycles on $P_{4} \times P_{10}$ according to the above weight. It turns out to be
$$
\matrix{
w_{1}^{9} w_{2}^{5} w_{3}^{5}+36 w_{1}^{9} w_{2}^{4} w_{3}^{6}+126 w_{1}^{9} w_{2}^{3} w_{3}^{7}+84 w_{1}^{9} w_{2}^{2} w_{3}^{8}+9 w_{1}^{9} w_{2} w_{3}^{9}+4 w_{1}^{8} w_{2}^{6} w_{3}^{5} 
+7 w_{1}^{8} w_{2}^{5} w_{3}^{6}+178 w_{1}^{8} w_{2}^{4} w_{3}^{7} \cr
+259 w_{1}^{8} w_{2}^{3} w_{3}^{8}+84 w_{1}^{8} w_{2}^{2} w_{3}^{9}+6 w_{1}^{7} w_{2}^{7} w_{3}^{5}  +42 w_{1}^{7} w_{2}^{6} w_{3}^{6}+137 w_{1}^{7} w_{2}^{5} w_{3}^{7}+178 w_{1}^{7} w_{2}^{4} w_{3}^{8} \cr
+126 w_{1}^{7} w_{2}^{3} w_{3}^{9}+4 w_{1}^{6} w_{2}^{8} w_{3}^{5} +15 w_{1}^{6} w_{2}^{7} w_{3}^{6}+42 w_{1}^{6} w_{2}^{6} w_{3}^{7}+67 w_{1}^{6} w_{2}^{5} w_{3}^{8}+36 w_{1}^{6} w_{2}^{4} w_{3}^{9}
+w_{1}^{5} w_{2}^{9} w_{3}^{5} \cr +4 w_{1}^{5} w_{2}^{8} w_{3}^{6}+6 w_{1}^{5} w_{2}^{7} w_{3}^{7}+4 w_{1}^{5} w_{2}^{6} w_{3}^{8}+w_{1}^{5} w_{2}^{5} w_{3}^{9}}
\quad .
$$

In particular, there are exactly $126$ Hamiltonian cycles of $P_4 \times P_{10}$ whose corresponding Stoyan-Strehl  $3 \times 9$ matrices have $9$ ones in the first row, $3$ ones in the second row,
and $7$ ones in the third row. (Note that the number of ones in the top row is the same as the number of edges visited at the upper boundary of the $P_m \times P_n$ grid, and the
number of ones in the bottom row is the number of edges visited at the bottom boundary of the grid).

Another example:  How do you find the {\bf exact} number of Hamiltonian cycles in $P_4 \times P_{100}$ whose corresponding $3 \times 99$ Stoyan-Strehl binary matrices have
$90$ ones in the first row, $31$ ones in the second row, and $78$ ones in the third row? Just type:

{\tt coeff(coeff(coeff(coeff(taylor(f,z=0,101),z,100),w[1],90),w[2],31),w[3],78);}

and you would get, {\it right away},
$$
1113455025360859674900898483836789708 \quad .
$$

As $m$ grows, these rational functions in $z$ and in $w_1, \dots, w_{m-1}$ grow more complicated (and take a long time to compute), but if you want
to focus on the number of ones in a specific row, the third argument of {\tt GFnzG} can be a list of size $m-1$ where the rows that you are not interested in are replaced by $1$. For example

{\tt coeff(expand(coeff(taylor(GFnzG(6,z,[w,1,1,1,1]),z=0,101),z,100)),w,80);} \quad ,

gives the number of Hamiltonian cycles on $P_6 \times P_{100}$ where $80$ edges (out of the total $99$) on the top boundary are visited. For the record, this number happens to be
$$
5769998174321676578317324842520250953447414723592327870562345553858388042 \quad .
$$

{\bf Expectation and Variance and their Asymptotics}

By taking derivatives with respect to the various $w_i$ (or if you focus on the number of ones in a specific row, say the first, w.r.t. $w$), and then
plugging-in  $w=1$, we get a rational function in $z$ whose coefficients is the sum of the statistic. We also have procedures to compute the asymptotics
of the actual enumerating sequence, and the implied asymptotics for the expectation and variance. Type {\tt ezraA();} for a list of the available
procedures.

{\bf Random Generation}

Using the {\it Wilf methodology} [W] [NW] for random generation of combinatorial objects, and once we had the infra-structure for enumeration, we wrote
a procedure {\tt RandHP(n,m)}, that outputs, {\it uniformly-at-random}, a Hamiltonian cycle on $P_m \times P_n$. For example

{\tt C:=RandHP(8,10); PlotC(C);}

To see many examples, look in this directory:

{\tt https://sites.math.rutgers.edu/\~{}zeilberg/tokhniot/volkerPics} \quad.

\vfill\eject

{\bf Output Files}

$\bullet$ For the generating functions for enumerating Hamiltonian Cycles in $P_m \times P_n$, for $2 \leq m \leq 9$, as well, for each, the first $50$ terms of the
enumerating sequence (all in the OEIS), see here:

{\tt https://sites.math.rutgers.edu/\~{}zeilberg/tokhniot/oVolker1.txt} \quad .

(Note that these are also included in the actual package, by typing {\tt GFnzPC(m,z);}) for $2 \leq m \leq 9$.)

$\bullet$ For the very complicated case of $m=10$, see here:

{\tt https://sites.math.rutgers.edu/\~{}zeilberg/tokhniot/oVolker3.txt} \quad .

(Note that it is also included in the actual package, by typing {\tt GFnzPC(10,z);}) \quad .

$\bullet$ For the generating function of the {\bf full} weight-enumerators, for $2 \leq m \leq 5$, see here:

{\tt https://sites.math.rutgers.edu/\~{}zeilberg/tokhniot/oVolker2.txt} \quad .

$\bullet$ If you want to see empirical confirmation for the Wilf-inspired procedure for random generation, compared to the
theoretical values of the expectation and variance see here:

{\tt https://sites.math.rutgers.edu/\~{}zeilberg/tokhniot/oVolker4.txt} \quad .

$\bullet$ If you want to see precise asymptotics for the average number of ones in various rows of the corresponding Stoyan-Strehl binary matrices, see here:

{\tt https://sites.math.rutgers.edu/\~{}zeilberg/tokhniot/oVolker5.txt} \quad.

$\bullet$ If you want to see precise asymptotics for the correlation between the number of edges visited on the top border and the bottom border, for
random Hamiltonian cycles in $P_m \times P_n$ for $3 \leq m \leq 7$ and $n \rightarrow \infty$,  see here:

{\tt https://sites.math.rutgers.edu/\~{}zeilberg/tokhniot/oVolker6.txt} \quad.

{\bf Conclusion: The Method is the Message}

With all due respect to the {\it actual content} that we outputted, we believe that far more important is the {\it method}, and
the present project is a nice
{\it case study} in the efficient and integrated  and streamlined use of symbolic computation. We combined, {\it automatically} and {\it seamlessly},
diverse tools and methods, in our {\bf favorite} computer algebra system, {\tt Maple}. It also brings to the attention of the enumerative
combinatorics community the wonderful, $30$-year-old  gem, of  Robert Stoyan (who had lots of potential, but apparently deserted academia) and
the {\it birthday boy}, Volker Strehl, an outstanding enumerator!

\vfill\eject

{\bf References}

[F] F.  Faase, {\it Counting Hamiltonian cycles in product graphs},\hfill\break
{\tt http://www.iwriteiam.nl/counting.html} \quad.

[Kr] Germain Kreweras, {\it D\'enombrement des cycles hamiltoniens dans un ractangle quadrill\'e}, European J. of Combinatorics {\bf 13} (1992), 473-467.

[KJ] A. Kloczkowski and R. L. Jernigan, {\it Transfer matrix method for enumeration and generation of compact self-avoiding walks. I. Square lattices}, The Journal of Chemical Physics {\bf 109} (1998), 5134.

[Kw] Harris Kwong, {\it Enumeration of Hamiltonian cycles in $P_4 \times P_n$ and $P_5 \times P_n$}, Ars Combinatorica {\bf 33} (1992), 87-96.

[KR]  Harris Y.H. Kwong and  D. G. Rogers, {\it A matrix method for counting Hamiltonian cycles on grid graphs},
European J. Combin. {\bf 15} (1994), 277-283.

[NW] Albert Nijenhuis and Herbert S. Wilf, “Combinatorial Algorithms”, Academic Press. First
edition: 1975. Second edition: 1978.\hfill\break
{\tt https://www2.math.upenn.edu/\~{}wilf/website/CombinatorialAlgorithms.pdf} \quad .

[Sl] Neil Sloane, The On-Line encyclopedia of integer sequences, {\tt https://oeis.org/} . Sequences A79, A6864, A6865, A145401, 145416, A145418, A160149, and  A180504.

[St] Richard P. Stanley, {\it Enumerative Combinatorics}, Volume 1, Wadsworth \& Brooks/Cole, (first edition), 1986.

[SS] Robert Stoyan and Volker Strehl, {\it Enumeration of Hamiltonian Circuits in Rectangular Grids}, S\'eminaire Lotharingien de Combinatoire, B34f (1995), 21pp. \hfill\break
{\tt https://www.mat.univie.ac.at/\~{}slc/wpapers/s34erlangen.html} \quad. \hfill\break
(Also published here: J. Combin. Math. Combin. Comput. {\bf 21} (1996), 109-127.)

[W] Herbert S. Wilf, A unified setting for sequencing, ranking, and selection algorithms for combinatorial objects, Advances in Math 24 (1977) , 281-291. \hfill\break
{\tt https://www2.math.upenn.edu/\~{}wilf/website/Unified\%20setting.pdf} \quad .
\bigskip
\hrule
\bigskip
Pablo Blanco and Doron Zeilberger, Department of Mathematics, Rutgers University (New Brunswick), Hill Center-Busch Campus, 110 Frelinghuysen
Rd., Piscataway, NJ 08854-8019, USA. \hfill\break
Emails: {\tt  pablancoh at aol dot com}  {\tt DoronZeil at gmail  dot com}   \quad .
\bigskip
{\bf Mar. 24, 2026} 

\end